\documentclass[12pt]{article}
\usepackage[margin=1.25in]{geometry}
\usepackage{amssymb,amsmath,amsthm,graphicx,enumitem,systeme,array,enumitem,commath,mathtools,multicol,url}
\usepackage{setspace,blindtext}
\usepackage{varioref}
\usepackage[colorlinks=true,linkcolor=blue]{hyperref}
\usepackage{cleveref}
\usepackage{blindtext}
\usepackage{fancyhdr}
\usepackage{pgfplots}
\usepackage[most]{tcolorbox}

\allowdisplaybreaks 

\usepackage{titlesec}
\titleformat*{\section}{\large\bfseries}

\newtheorem{theorem}{Theorem}[section]

\newtheorem{lemma}[theorem]{Lemma}
\newtheorem{claim}[theorem]{Claim}
\newtheorem{conjecture}[theorem]{Conjecture}
\newtheorem{remark}[theorem]{Remark}
\newtheorem*{remark*}{Remark}

\newtheorem{definition}[theorem]{Definition}
\newtheorem{example}[theorem]{Example}

\begin{document}


\begin{center}
    \Large{{\textbf{Variants of the Damascus inequality}}}
    \\[25pt]
\end{center}
\begin{center}
    \textsc{Chanatip Sujsuntinukul} and \textsc{Christophe Chesneau}
\end{center}

\begin{abstract}
In 2016, Dannan and Sitnik established the notable Damascus inequality, which features a symmetric structure under a multiplicative constraint. In this study, we consider the natural generalisation of this inequality by characterising all positive integers $m$ and $n$ such that the inequality
\[\sum_{j=1}^m\frac{x_j^n-1}{x_{j}^{n+1}+1}\leqslant 0\]
holds for any positive real numbers $x_1, \ldots, x_m$ with $\prod_{j=1}^mx_j=1$.  
Our approach relies on the theories of GA-convexity and Sturm's sequence.
For the cases where the inequality fails, we also investigate the topological properties of the set of non-solutions. 
\end{abstract}




\section{Introduction}

Inequalities form the backbone of mathematical analysis and optimisation, with symmetry and convexity often determining both their structure and sharpness. Classical examples illuminate these principles: the arithmetic-geometric mean inequality (AM-GM) states that for nonnegative real numbers
$x_1, \ldots, x_n$, 
\[\frac{1}{n}\sum_{j=1}^nx_j\geqslant\left(\prod_{j=1}^nx_j\right)^{1/n},\]
with equality if and only if $x_1=\cdots=x_n$ \cite{BEF, HGHL}.
Similarly fundamental results—such as the Cauchy-Schwarz inequality \cite{SJM}, Jensen's inequality for convex functions \cite{JJLWV}, and Chebyshev's sum inequality for ordered sequences \cite{MDS}—share common threads of symmetry, convexity, and majorisation that continue to guide modern developments in inequality theory.

Within this rich landscape, the Damascus inequality emerged in 2016 as a surprising result. Proposed by Dannan from Damascus, Syria, it states that for positive real numbers $x$, $y$, and $z$ satisfying $xyz=1$, 
\begin{align*}
    \frac{x-1}{x^2-x+1}+\frac{y-1}{y^2-y+1}+\frac{z-1}{z^2-z+1}\leqslant 0,
\end{align*}
with equality if and only if $(x, y, z)=(1, 1, 1)$. Despite its elementary appearance, this inequality resisted straightforward verification and sparked considerable interest. Dannan and Sitnik \cite{FMD} eventually established it using Lagrange multipliers, cyclic inequalities, and geometric methods, while Kalinin and Pankratova \cite{SIK} provided an alternative proof employing geometric means, logarithmic means, and GA-convex functions.

Building on these foundations, Kalinin and Pankratova posed a natural generalisation: for which positive integers $n$ does the inequality\[ \frac{x^n-1}{x^{n+1}+1}+\frac{y^n-1}{y^{n+1}+1}+\frac{z^n-1}{z^{n+1}+1}\leqslant 0\]hold for positive real numbers $x$, $y$, $z$ with $xyz = 1$? Moreover, when does the inequality become an equality? They explicitly resolved only the cases $n=1, 2$, leaving the general question open. 

Questions of exponent extension as such appear widely in inequality theory, often revealing deeper structural properties and delimiting the scope of particular techniques. A systematic characterisation of all valid exponents not only completes the theoretical picture but also provides practical guidance for applications in optimisation and analysis. For our setting, the difficulty is reflected by the fact that the function
\[
f(x) := \frac{1 - e^{nx}}{e^{(n+1)x} + 1}
\]
is not half convex on $\mathbf{R}$ (for more information, read \cite{ZP}).

In this paper, we provide a complete answer to the question posed by Kalinin and Pankratova concerning this extension. In fact, we introduce and fully resolve a related question in a more general setting. Particularly, we characterise positive integers $m$ and $n$ (provided in Table \ref{mainresult}) such that the inequality 
\begin{equation}\label{varineq}
    S_{n}^m(x_1, \ldots, x_m):=\sum_{j=1}^m\frac{x_j^n-1}{x_{j}^{n+1}+1}\leqslant 0
\end{equation}
holds on the $(m-1)$-dimensional smooth constraint hypersurface
\[\mathcal{H}_m:=\left\{(x_1, \ldots, x_m)\in(0, +\infty)^m : \prod_{j=1}^mx_j=1\right\}.\]

\begin{table}[]
    \centering
     \begin{tabular}{c|c}
      $m$ & $n$ \\
      \hline
       $1$  & every positive integer \\
       $2$  & every positive integer \\
       $3$ & $1$, $2$, $3$ \\
       $4$ & $1$ \\
       $5$ & $1$
    \end{tabular}
    \caption{Main result}
    \label{mainresult}
\end{table}

Furthermore, for $m=3$ we study the structure of the violation set under the original constraint, proving that it is bounded and separated from certain planes.

\section{Cases $m=1, 2$}\label{subsec2.1}

The case $m=1$ is straightforward (true). When $m=2$, we can show that \eqref{varineq} holds for all values of $n$ using a basic calculus argument. We demonstrate this below.

\begin{theorem}
For any $n\in\mathbf{N}$,
$S_{n}^2\leqslant 0$ on $\mathcal{H}_2$ and $(S_n^2)^{-1}(\{0\})=\{(1, 1)\}$.
\end{theorem}
\noindent {\bf Proof.} 
By some simplification and sign argument, it is equivalent to consider the inequality
    \begin{align*}
        (x^n-1)(y^{n+1}+1)+(y^n-1)(x^{n+1}+1)\leqslant 0,
    \end{align*}
    which is the same as
    \[x^ny^{n+1}+x^n-y^{n+1}+x^{n+1}y^n+y^n-x^{n+1}-2\leqslant 0.\]
    Substituting the constraint $xy=1$ gives
    \[\frac{1}{x}+x^n-\frac{1}{x^{n+1}}+x+\frac{1}{x^n}-x^{n+1}-2\leqslant 0.\]
    For convenience, define $\psi : (0, +\infty)\to\mathbf{R}$ by the formula
    \[\psi(x)=\frac{1}{x}+x^n-\frac{1}{x^{n+1}}+x+\frac{1}{x^n}-x^{n+1}-2.\]
    Then we have 
    \begin{align*}   \psi'(x) &=-\frac{1}{x^2}+nx^{n-1}+\frac{n+1}{x^{n+2}}+1-\frac{n}{x^{n+1}}-(n+1)x^n\\     &= \left(\frac{1}{x^{n+2}}-\frac{1}{x^2}\right)+n(x^{n-1}-x^n)+n\left(\frac{1}{x^{n+2}}-\frac{1}{x^{n+1}}\right)+(1-x^n).\end{align*}
    There are three parts to consider.
    \begin{itemize}
        \item For $x<1$, we note the following inequalities.
    \begin{align}\label{someineq}
        \frac{1}{x^{n+2}}>\frac{1}{x^2},\quad nx^{n-1}>nx^n, \quad \frac{n}{x^{n+2}}>\frac{n}{x^{n+1}}, \quad 1>x^n.
    \end{align}
    This shows $\psi'(x)>0$.
    \item Next, we have $\psi'(1)=0$ and $\psi(1)=0$.
    \item Lastly, if $x>1$, we see that inequalities in \eqref{someineq} hold when the signs are swapped from $>$ to $<$.
    This means $\psi'(x)<0$.
    \end{itemize}
    Thus by the first derivative test, the function 
$\psi$ has maximum value $0$ uniquely at $x=1$. In this case $(x, y)=(1, 1)$.
    So the proof is complete.\hfill $\square$

\section{Case $m=3$}\label{subsec2.2}

For the case $m=3$, we know in the previous literature that \eqref{varineq} holds when $n=1, 2$. Here,
we first find that the inequality fails for $n\geqslant 4$.

\begin{theorem}\label{m=3n>=4}
Let $n\in\mathbf{N}$ such that $n\geqslant 4$. Then there exists $(x, y, z)\in\mathcal{H}_3$ such that $S_n^3(x, y, z)>0$.
\end{theorem}
\noindent {\bf Proof.} 
    For $n=4$, take $(x, y, z)=(31/20, 43/100, 2000/1333)$ as an example. It is easy to check that the constraint $xyz=1$ is met. Then 
    \[S_n^3=\frac{354458009159794612949999}{481099388060786340236540521}>0.\]

Next, for $n\geqslant 5$, we claim that $(x, y, z)=(1/2, 3/2, 4/3)$ is an example for $S_n^3>0$. Indeed, we can check that the condition $xyz=1$ is satisfied. Equivalently, we have to show that
\begin{align}\label{eqn>=4_1}
  \frac{\left(\frac{1}{2}\right)^{n}-1}{\left(\frac{1}{2}\right)^{n+1}+1}+\frac{\left(\frac{3}{2}\right)^{n}-1}{\left(\frac{3}{2}\right)^{n+1}+1}+\frac{\left(\frac{4}{3}\right)^{n}-1}{\left(\frac{4}{3}\right)^{n+1}+1} >0
\end{align}
for $n\geqslant 5$. First, by adding terms in \eqref{eqn>=4_1} and using the fact that the denominator is always positive, we need to show that
\begin{align*}   
\left(\left(\frac{1}{2}\right)^n-1\right)\left(\left(\frac{3}{2}\right)^{n+1}+1\right)\left(\left(\frac{4}{3}\right)^{n+1}+1\right)\\+\left(\left(\frac{3}{2}\right)^n-1\right)\left(\left(\frac{1}{2}\right)^{n+1}+1\right)\left(\left(\frac{4}{3}\right)^{n+1}+1\right)\\+\left(\left(\frac{4}{3}\right)^n-1\right)\left(\left(\frac{1}{2}\right)^{n+1}+1\right)\left(\left(\frac{3}{2}\right)^{n+1}+1\right)>0.
\end{align*}
Multiplying by $2^{2n+2}\cdot3^{n+1}$, the inequality becomes
\begin{align*}
   2(1-2^n)(3^{n+1}+2^{n+1})(4^{n+1}+3^{n+1})+2(3^n-2^n)(1+2^{n+1})(4^{n+1}+3^{n+1})\\
   +3(4^n-3^n)(1+2^{n+1})(3^{n+1}+2^{n+1}) >0.
\end{align*}
By expanding the left-hand side expression and performing some algebraic manipulation, we need to prove that
\begin{align}\label{eqn>=4}
\begin{split}
   7\cdot 2^{3n+1}-5\cdot 2^{4n+2}+5\cdot 2^{2n}\cdot 3^n+5\cdot 2^{3n+1}\cdot 3^n\\+5\cdot 3^{2n+1}-2^{n+3}\cdot 3^{2n+1}>0.
\end{split}
\end{align}
Now, we observe that
\begin{align*}
  \frac{2^{n+3}\cdot 3^{2n+1}+5\cdot 2^{4n+2}}{5\cdot 2^{3n+1}\cdot 3^n}=\frac{12}{5}\left(\frac{3}{4}\right)^n+2\left(\frac{2}{3}\right)^n<1
\end{align*}
for $n\geqslant 5$. This can easily be verified by checking that the inequality is true for $n=5$ and that $(3/4)^n$ and $(2/3)^n$ are decreasing sequences. Thus  
\[5\cdot 2^{3n+1}\cdot 3^n-2^{n+3}\cdot 3^{2n+1}-5\cdot 2^{4n+2}> 0.\] Hence \eqref{eqn>=4} is true and $S_n^3>0$ for $n\geqslant 5$ and this triple.\hfill $\square$

\hfill

Let us now consider when $n=3$. We will be showing that \eqref{varineq} holds. 

\begin{theorem}\label{finally}
We have $S_3^3\leqslant 0$ on $\mathcal{H}_3$ and $(S_3^3)^{-1}(\{0\})=\{(1, 1, 1)\}$.
\end{theorem}

Whilst one might consider applying the Lagrange multiplier method to resolve this inequality directly, such an approach proves exceedingly laborious owing to the degree of the variables involved. Hence the technique we will be using is the same as in \cite{SIK}. First, we need to introduce several useful terminologies and lemmas related to the geometric-arithmetic convexity. They can be found in \cite{CPN, CPN2}.

\begin{definition}
    Let $I\subseteq(0, +\infty)$ be an interval, and let $f : I\to\mathbf{R}$. We say that $f$ is \textit{GA-convex} on $I$ if for any $\alpha, \beta\in I$ and any $\lambda\in[0, 1]$, we have
    \begin{equation}\label{GA}
        f(\alpha^{\lambda}\beta^{1-\lambda})\leqslant \lambda f(\alpha)+(1-\lambda)f(\beta).
    \end{equation}
    If \eqref{GA} is reversed, then we say that $f$ is \textit{GA-concave} on $I$. 
\end{definition}

\begin{remark}\label{rem1}
     On $I$, $f$ is GA-convex if and only if $f(e^x)$ is convex.
\end{remark}

\begin{remark}
One can prove (by induction) that if $f$ is GA-convex on $I$, then for any $\alpha_1, \ldots, \alpha_n\in I$ and any $\lambda_1, \ldots, \lambda_n\in[0, 1]$ with $\sum_{j=1}^n\lambda_j=1$,
\[f\left(\prod_{j=1}^n \alpha_j^{\lambda_j}\right)\leqslant\sum_{j=1}^n\lambda_jf(\alpha_j).\]
If $f$ is GA-concave, then the inequality is reversed. This fact will be useful when we consider a large $m$ case.
\end{remark}

\begin{lemma}\label{GAconlem}
    Let $I\subseteq(0, +\infty)$ be an interval, and let $f : I\to\mathbf{R}$ be a continuous function  such that $f$ is twice differentiable on $I^{\circ}$. Define 
    \[\Delta_f(x)=f'(x)+xf''(x)\]
    for $x\in I^{\circ}$. If $\Delta_f(x)\geqslant 0$ for all $x\in I$, then $f$ is GA-convex on $I$. On the other hand, if $\Delta_f(x)\leqslant 0$, then $f$ is GA-concave on $I$.
\end{lemma}

Let us now prove the result we want.

\noindent {\bf Proof of Theorem \ref{finally}.} 
    First of all, define the function $f : (0, +\infty)\to\mathbf{R}$ by the formula
\[f(x)=\frac{x^3-1}{x^4+1}-\frac{3}{2}\ln x.\]
We compute
\[f'(x)=-\frac{3x^8+2x^7-2x^4-6x^3+3}{2x(x^4+1)^2}\]
to find the critical points of $f$. Next, we may build up  
Sturm's sequence of the function $-2x(x^4+1)^2f'(x)$, using a computational software\footnote{One may access Maple or the website \url{https://planetcalc.com/7719/}. For more information regarding Sturm's sequence, we refer the readers to \cite{KEA}.}. We demonstrate this below.
\begin{align*}
    f_0(x)&=-2x(x^4+1)^2f'(x)\\
    f_1(x)&=24x^7+14x^6-8x^3-18x^2\\
    f_2(x)&=\frac{7}{48}x^6+x^4+\frac{11}{3}x^3-\frac{3}{16}x^2-3\\
    f_3(x)&=\frac{1152}{7}x^5+\frac{4896}{7}x^4+\frac{2304}{7}x^3-\frac{3456}{7}x-288\\
    f_4(x)&=-\frac{2567}{768}x^4-\frac{157}{32}x^3-\frac{1}{4}x^2+\frac{77}{48}x+\frac{3137}{768}\\
    f_5(x)&=\frac{2340864000}{6589489}x^3-\frac{294801408}{6589489}x^2+\frac{480079872}{6589489}x-\frac{1789231104}{6589489}\\
    f_6(x)&=\frac{728975399603}{3096768000000}x^2-\frac{146490929959}{1032256000000}x-\frac{1653961739}{129032000000}\\
    f_7(x)&=-\frac{1424761021440000000}{7331305593807371}x+\frac{113118741504000000}{431253270223963}\\
    f_8(x)&=-\frac{31931834754991359271}{142253483234400000000}
\end{align*}
Then we construct the table of signs of this finite sequence in Table \ref{tablesturm}.
\begin{table}[]
    \centering
    \begin{tabular}{c|c|c}
             & $0$ & $+\infty$  \\
             \hline
             $f_0$ & $+$ & $+$ \\
             $f_1$ & $0$ & $+$ \\
             $f_2$ & $-$ & $+$ \\
             $f_3$ & $-$ & $+$ \\
             $f_4$ & $+$ & $-$ \\
             $f_5$ & $-$ & $+$ \\
             $f_6$ & $-$ & $+$ \\
             $f_7$ & $+$ & $-$ \\
             $f_8$ & $-$ & $-$ \\
             \hline
           $V$  & $5$ & $3$
        \end{tabular}
    \caption{Sturm's sequence sign table of $-2x(x^4+1)^2f'(x)$}
    \label{tablesturm}
\end{table}

Thus we yield $V(0)-V(+\infty)=5-3=2$. 
This shows that $f'$ has precisely two positive real roots: one root is precisely at $x=1$, another, denoted by $\alpha$, lies in $(4/5, 9/10)$ by the intermediate value theorem.
Note that \(f(1)=0\) and \(f(4/5), f(\alpha)<0\). Hence, since \(f\) is continuous, satisfies \(\lim_{x\to+\infty}f(x)=-\infty\), and $-2x(x^4+1)^2f'(x)$ has a single root after \(\alpha\), we deduce that \(f(x)\leqslant 0\) for all \(x\in[4/5, +\infty)\).

In the next step, we aim to find an upper bound for $f(x)+f(y)+f(z)$. First, if $x, y, z\in[4/5, +\infty)$ and $xyz=1$, then \[f(x)+f(y)+f(z)\leqslant 0,\] which proves our claim (note that $\ln x+\ln y+\ln z=\ln xyz=0$), and the equality is achieved only if $f(x)=f(y)=f(z)=0$, i.e., $x=y=z=1$.

Second, if some value among $x, y, z$ does not belong to the interval $[4/5, +\infty)$, then there are two cases to consider (up to symmetry). For convenience, we define a function $h : [0, +\infty)\to\mathbf{R}$ by the formula
\[h(x)=\frac{x^3-1}{x^4+1}.\]
Observe that $h(x)+h(y)+h(z)=f(x)+f(y)+f(z)$ for $x, y, z>0$ with $xyz=1$ by the earlier discussion. We compute
\[h'(x)=-\frac{x^2(x^4-4x-3)}{(x^4+1)^2},\]
which tells us that $h$ is increasing on $(0, 1]$. So $h(x)<h(4/5)=-305/881$ for $0<x<4/5$.
Also, we claim that
$\max\{h(x) : x\in[0, +\infty)\}\in (4203/10000, 1051/2500)$
and it is attained at some  $x\in (223/125, 357/200)$. To roughly show this, we first note that $h(x)\geqslant 0$ when $x\geqslant 1$, and $\lim_{x\to+\infty}h(x)=0$. Since $h$ is continuous, we can argue using the extreme value theorem that $h$ has a global maximum. Next, we find critical points of $h$ in $[0, +\infty)$. One can verify that $h'$ has two nonnegative real roots (e.g. by observing its monotonicity), one at $x=0$, and another at some $x\in (223/125, 357/200)$ by the intermediate value theorem. By comparing the images of these points, we prove the claim. Now, let us turn to the case analysis.

\textit{Case 1:} $x, y\in(0, 4/5)$ and $z\in[4/5, +\infty)$. 
We have
\[h(x)+h(y)+h(z)< -\frac{2(305)}{881}+\frac{1051}{2500}<0.\]

\textit{Case 2:} $x\in (0, 4/5)$, $y, z\in[4/5, +\infty)$, $y\leqslant z$. If $x\leqslant 2/5$, then we are done because 
\begin{align*}
    h(x)+h(y)+h(z)< h\left(\frac{2}{5}\right)+2\left(\frac{1051}{2500}\right)<0.
\end{align*}
If $y\leqslant 9/10$, then we are also done since \begin{align*}
    h(x)+h(y)+h(z) &< h\left(\frac{4}{5}\right)+h\left(\frac{9}{10}\right)+\frac{1051}{2500}\\
    &=-\frac{305}{881}-\frac{2710}{16561}+\frac{1051}{2500}<0.
\end{align*}
Therefore we are only left to consider the case  $x\in(2/5, 4/5)$ and $y, z\in(9/10, +\infty)$ with $y\leqslant z$. However, due to the constraint $xyz=1$, we are forced to have $y, z\in(9/10, 14/5)$.
Let us analyse the GA-convexity of $h(y)$. We see that
\begin{align*}
    \Delta_h(y) &=h'(y)+yh''(y)\\
    &=-\frac{y^2(y^4-4y-3)}{(y^4+1)^2}+\frac{2y^2(y^8-10y^5-12y^4+6y+3)}{(y^4+1)^3}\\
    &=\frac{y^2(y^8-16y^5-22y^4+16y+9)}{(y^4+1)^3}.
\end{align*}
We may argue that $\Delta_h(y)<0$ for $y\in(9/10, 14/5)$ by constructing Sturm's sequence of $(y^4+1)^3\Delta_h(y)$ and seeing that it has no real roots in this interval.
Combining this fact, the continuity of $\Delta_h(y)$, and $\Delta_h(1)<0$, we deduce that $h$ is GA-concave on $(9/10, 14/5)$ according to Lemma \ref{GAconlem}. This implies
\begin{align*}
    h(y)+h(z)\leqslant 2h(\sqrt{yz})=2h\left(\frac{1}{\sqrt{x}}\right)
    =\frac{2\sqrt{x}(1-x\sqrt{x})}{x^2+1}.
\end{align*}
Therefore in the condition $x\in(2/5, 4/5)$, $9/10<y\leqslant z<14/5$, we have
\begin{align*}
    h(x)+h(y)+h(z) &\leqslant\frac{x^3-1}{x^4+1}+\frac{2\sqrt{x}(1-x\sqrt{x})}{x^2+1} \\
    &=\frac{(x^3-1)(x^2+1)+2\sqrt{x}(1-x\sqrt{x})(x^4+1)}{(x^2+1)(x^4+1)} \\
    &=-\frac{(\sqrt{x}-1)^2(x+\sqrt{x}+1)\ell(x)}{(x^2+1)(x^4+1)},
\end{align*}
where $\ell(x):=2x^4+2x^3\sqrt{x}+x^3+x^2\sqrt{x}+x^2-x-\sqrt{x}+1$. One can verify that $\ell$ is positive on $(2/5, 4/5)$ (e.g. by making a substitution $t=\sqrt{x}$, then constructing  Sturm's sequence). 
Thus $h(x)+h(y)+h(z)<0$ and the proof is complete.\hfill $\square$


Now, although $(1, 1, 1)$ is not a global maximum of $S_n^3$
under $\mathcal{H}_3$ for $n\geqslant 4$, it is a (strict) local maximum. This can be shown by applying the second derivative test to the function $S_n^3(x, y, 1/xy)$. A direct computation shows
\begin{align*}
    \partial_x S_n^3 &=\frac{n}{x(x^{-n}+x)}-\frac{(n+1)(x^n-1)x^n}{(x^{n+1}+1)^2}-\frac{ny}{(xy)^{n+1}+1} -\frac{(n+1)y((xy)^n-1)}{((xy)^{n+1}+1)^2},\\
    \partial_y S_n^3 &=\frac{n}{y(y^{-n}+y)}-\frac{(n+1)(y^n-1)y^n}{(y^{n+1}+1)^2}-\frac{nx}{(xy)^{n+1}+1} -\frac{(n+1)x((xy)^n-1)}{((xy)^{n+1}+1)^2}.
\end{align*}
This gives
$\nabla S_n^3(1, 1, 1)=\mathbf{0}$. A further calculation shows
\[HS_n^3(1, 1, 1)=\begin{pmatrix}
    -n & -\frac{n}{2}\\
    -\frac{n}{2} & -n
\end{pmatrix}.\]
The principal minors of $HS_n^3(1, 1, 1)$ are 
\[d_1=-n<0\quad \text{and}\quad d_2=n^2-\left(\frac{n}{2}\right)^2=\frac{3n^2}{4}>0.\]
Hence by Sylvester's criterion, $HS_n^3(1, 1, 1)$ is negative definite. Thus the claim follows.


\section{Cases $m\geqslant 4$}

Next,  for \(m\geqslant 4\), the analysis becomes more intricate due to the number of variables; nevertheless, a large number of these cases can be ruled out, as illustrated below.

\begin{theorem}\label{counterexample}
Let $m, n\in\mathbf{N}$. If $m\geqslant 4$ and $n\geqslant 2$ or if $m\geqslant 6$ and $n\geqslant 1$, then there exists $(x_1, \ldots, x_m)\in\mathcal{H}_m$ such that $S_n^m(x_1, \ldots, x_m)>0$.
\end{theorem}
\noindent {\bf Proof.} 
    Our choice of counterexample is inspired by \cite{FMD}. For the first case, we consider $(x_1, \ldots, x_m)=(2, 2, \ldots, 2, 1/2^{m-1})$. Formally, every $x_j$ except $x_m$ is $2$. Trivially, the multiplicative constraint is met. Now, we have
    \begin{align*}
        \sum_{j=1}^m\frac{x_j^n-1}{x_{j}^{n+1}+1}&=\sum_{j=1}^{m-1}\frac{2^n-1}{2^{n+1}+1}+\frac{(1/2^{m-1})^n-1}{(1/2^{m-1})^{n+1}+1} \\
        &=\frac{(m-1)(2^n-1)}{2^{n+1}+1}+\frac{2^{n(1-m)}-1}{2^{(n+1)(1-m)}+1} \\
        &=\frac{(m-1)(2^n-1)(2^{(n+1)(1-m)}+1)+(2^{n(1-m)}-1)(2^{n+1}+1)}{(2^{n+1}+1)(2^{(n+1)(1-m)}+1)}.
    \end{align*}
    Let us consider the numerator of the expression above as the denominator is always positive. By simplification, it equals
    \begin{align*}
      f(n, m)&:= m2^n-m2^{(1-m)(n+1)}+m2^{(1-m)(n+1)+n}+2^{(1-m)n}-2^n\\
      &\quad+2^{(1-m)(n+1)}+2^{(1-m)n+n+1}-2^{(1-m)(n+1)+n}-2^{n+1}-m.
    \end{align*}
   We see that 
   \begin{align*}
       m2^n-2^n-2^{n+1}-m &=m(2^n-1)-2^n-2^{n+1}\geqslant 4(2^{n}-1)-2^n-2^{n+1}>0
   \end{align*}
   for $m\geqslant 4$. Moreover, we have
   \begin{align*}
       &m2^{(1-m)(n+1)+n}-m2^{(1-m)(n+1)}-2^{(1-m)(n+1)+n} \\
       &\quad=(m-1)2^{(1-m)(n+1)+n}-m2^{(1-m)(n+1)} \\
       &\quad\geqslant (4m-4)2^{(1-m)(n+1)}-m2^{(1-m)(n+1)} \\
       &\quad= (3m-4)2^{(1-m)(n+1)}>0
   \end{align*}
   for $m\geqslant4$ and $n\geqslant 2$. Hence we conclude that $f(n, m)$ is positive on the given domain.
   
    For the second case, it suffices to only verify the cases $m\geqslant 6$ and $n=1$. We take the same counterexample. We have
    \begin{align*}
        \sum_{j=1}^m\frac{x_j-1}{x_{j}^{2}+1}&=\sum_{j=1}^{m-1}\frac{2-1}{2^{2}+1}+\frac{(1/2^{m-1})-1}{(1/2^{m-1})^{2}+1} \\
        &=\frac{m-1}{5}+\frac{2^{1-m}-1}{2^{2(1-m)}+1} \\
        &=\frac{(m-1)(2^{2(1-m)}+1)+5(2^{1-m}-1)}{5(2^{2(1-m)}+1)}.
    \end{align*}
    Since the denominator of the expression above is always positive, we consider the numerator:
    \begin{align*}
        g(m) &:=(m-1)(2^{2(1-m)}+1)+5(2^{1-m}-1) \\
        &=(m-1)2^{2(1-m)}+5(2^{1-m})+(m-6)>0
    \end{align*}
    since $m\geqslant 6$ and $(m-1)2^{2(1-m)}+5(2^{1-m})>0$.\hfill $\square$

\begin{remark}\label{rmk1}
Theorem \ref{counterexample} admits a more streamlined proof. It suffices to construct counterexamples only for the cases $(m, n)\in(\{4\}\times\mathbf{N}_{\ge 2})\cup\{(6, 1)\}$. The remaining cases then follow by setting the appropriate variables to $1$.
\end{remark}

We now consider the case $(m, n)=(5, 1)$ and establish its validity. The proof follows the same approach as that employed for the case $(m, n)=(3, 3)$; consequently, we present only a proof sketch, omitting certain routine computational details. In addition, as a corollary of this, \eqref{varineq} holds for $(m, n)=(4, 1)$.

\begin{theorem}\label{m=5n=1}
    We have $S_1^5\leqslant 0$ on $\mathcal{H}_5$ and $(S^5_1)^{-1}(\{0\})=\{(1, 1, 1, 1, 1)\}$.
\end{theorem}
\noindent {\bf Proof.} 
    We begin by defining the function $f : (0, +\infty)\to\mathbf{R}$ using the formula
    \[f(x)=\frac{x-1}{x^2+1}-\frac{1}{2}\ln x.\]
    We compute
    \[f'(x)=\frac{-x^4-2x^3+2x^2+2x-1}{2x(x^2+1)^2}\]
    to find the critical points of $f$. One can construct Sturm's sequence of $2x(x^2+1)^2f'(x)$ to show that it has two positive real roots, one at $x=\sqrt{2}-1$, another at $x=1$. We note that $f(117/500)\leqslant 0$. Combining these facts, we deduce that $f(x)\leqslant 0$ on $[117/500, +\infty)$ with equality at $x=1$.

    Next, let us try to estimate $f(s)+f(w)+f(x)+f(y)+f(z)$ from above. If $s, w, x, y, z\in[117/500, +\infty)$ and $swxyz=1$, then
    \[f(s)+f(w)+f(x)+f(y)+f(z)\leqslant 0.\]
    Since we know that $\ln swxyz=\ln 1=0$, this proves our claim, and the equality is attained precisely when $s=w=x=y=z=1$.

    Now, we consider when some value among $s, w, x, y, z$ is less than $117/500$. There are four cases to consider up to symmetry. Before we delve into each of them, we define a function $h : [0, +\infty)\to\mathbf{R}$ by the formula
    \[h(x)=\frac{x-1}{x^2+1}.\]
    Note that 
    \[h(s)+h(w)+h(x)+h(y)+h(z)=f(s)+f(w)+f(x)+f(y)+f(z)\] for any $s, w, x, y, z>0$ with $swxyz=1$. By computing
    \[h'(x)=\frac{-x^2+2x+1}{(x^2+1)^2}\]
    and using routine calculus argument, one can show that $h$ has a unique global maximum equals $1/\sqrt{2}-1/2$, attained at $x=1+\sqrt{2}$. Moreover, this function $h$ is increasing on $(0, 1+\sqrt{2})$ and decreasing on $(1+\sqrt{2}, +\infty)$.

     \textit{Case 1:} $s, w, x, y\in(0, 117/500)$ and $z\in[117/500, +\infty)$. Then
    \begin{align*}
        h(s)+h(w)+h(x)+h(y)+h(z)\leqslant 4h\left(\frac{117}{500}\right)+\left(\frac{1}{\sqrt{2}}-\frac{1}{2}\right)<0.
    \end{align*}
    
    \textit{Case 2:} $s, w, x\in(0, 117/500)$ and $y, z\in[117/500, +\infty)$. Then
    \begin{align*}
       h(s)+ h(w)+h(x)+h(y)+h(z)\leqslant 3h\left(\frac{117}{500}\right)+2\left(\frac{1}{\sqrt{2}}-\frac{1}{2}\right)<0.
    \end{align*}
    
    \textit{Case 3:} $s, w\in(0, 117/500)$ and $x, y, z\in[117/500, +\infty)$. Then 
    \begin{align*}
      h(s)+  h(w)+h(x)+h(y)+h(z)\leqslant 2h\left(\frac{117}{500}\right)+3\left(\frac{1}{\sqrt{2}}-\frac{1}{2}\right)<0.
    \end{align*}

    \textit{Case 4:} $s\in(0, 117/500)$ and $w, x, y, z\in[117/500, +\infty)$. Without loss of generality, suppose $w\leqslant x\leqslant y\leqslant z$.
    First, if $s\leqslant 3/20$, then we are done because
    \begin{align*}
        h(s)+  h(w)+h(x)+h(y)+h(z)\leqslant h\left(\frac{3}{20}\right)+4\left(\frac{1}{\sqrt{2}}-\frac{1}{2}\right)<0.
    \end{align*}
    If $w\leqslant 127/100$, then we are also done since
    \begin{align*}
       & h(s)+  h(w)+h(x)+h(y)+h(z)\\&
        \quad \leqslant h\left(\frac{117}{500}\right)+h\left(\frac{127}{100}\right)+3\left(\frac{1}{\sqrt{2}}-\frac{1}{2}\right)<0.
    \end{align*}
     Thus we consider when $s\in(3/20, 117/500)$ and $w, x, y, z\in(127/100, +\infty)$. However, due to the constraint $swxyz=1$, we are forced to have 
     \[w, x, y, z\in(127/100, 20000000/6145149).\]  Let us now compute
     \begin{align*}
         \Delta_h(w) &=-\frac{w^2-2w-1}{(w^2+1)^2}+\frac{2w(w^3-3w^2-3w+1)}{(w^2+1)^3} \\
         &=\frac{w^4-4w^3-6w^2+4w+1}{(w^2+1)^3}.
     \end{align*}
     One can check that $\Delta_h(w)<0$ for $w\in(127/100, 20000000/6145149)$ by
     using Sturm's sequence. Therefore $h$ is GA-concave in this range, and so we deduce that
     \begin{align*}
         h(w)+h(x)+h(y)+h(z)\leqslant 4h(\sqrt[4]{wxyz})=4h\left(\frac{1}{\sqrt[4]{s}}\right).
     \end{align*}
     So in the condition we are currently looking, we have
     \begin{align*}
         h(s)+h(w)+h(x)+h(y)+h(z)
         &\leqslant h(s)+4h\left(\frac{1}{\sqrt[4]{s}}\right) \\
         &= \frac{s-1}{s^2+1}+\frac{4(1/\sqrt[4]{s}-1)}{1/\sqrt{s}+1}.
     \end{align*}
     Finally, one can argue using calculus or Sturm's sequence that the above expression is negative for $s$ in the given domain. Thus the proof is complete.\hfill $\square$

\begin{remark}We note that the coefficient of the natural logarithmic function in the definitions of $f$ differs between the proofs of Theorem \ref{finally} and Theorem \ref{m=5n=1}. This difference arises from our requirement that $f$ possess a critical point at $x=1$.\end{remark}

\section{Non-solutions discussion}

Having characterised all pairs $(m, n)$ for which \eqref{varineq} holds, we now turn our attention to the non-solutions to this inequality. 
For convenience, we define the set
\[\mathcal{I}_n^m:=\left\{(x_1, \ldots, x_m)\in\mathcal{H}_m :  \sum_{j=1}^m\frac{x_j^n-1}{x_{j}^{n+1}+1}>0\right\}.\]
For instance, we know that $\mathcal{I}^2_n=\emptyset$ for all $n$, and $(1/2, 3/2, 4/3)\in\mathcal{I}^3_n$ for $n\geqslant 5$. A trivial observation is if 
$\mathcal{I}_n^m\neq\emptyset$, then 
$\mathcal{I}_n^m$ is a nonempty open subset of $\mathcal{H}_m$.
 Consequently, $\mathcal{I}_n^m$ has positive $(m-1)$-dimensional measure with respect to the induced Lebesgue measure on this surface.

Particularly for the case $m=3$, we find two topological properties for $\mathcal{I}_n^3$.

\begin{theorem}\label{3farfrom1}
    For each $n\in\mathbf{N}_{\ge 4}$, there exists $\delta>0$ such that $\mathcal{I}_n^3$ is disjoint from the set
    \[([1-\delta, 1+\delta]\times\mathbf{R}^2)\cup(\mathbf{R}\times[1-\delta, 1+\delta]\times\mathbf{R})\cup(\mathbf{R^2}\times[1-\delta, 1+\delta]).\]
\end{theorem}

This result shows that the set of non-solutions $\mathcal{I}_n^3$ stays at a positive distance from the loci where at least one of $x,y,z$ equals $1$. Building on this separation property, we can also derive the following consequence.

\begin{theorem}\label{3bound}
Fix $n\in\mathbf{N}_{\ge 4}$.
    The set $\mathcal{I}_n^3$ is bounded. Equivalently, there exists $c\in(0, 1]$ such that $\mathcal{I}_n^3\subseteq[c, 1/c]^3$.
\end{theorem}

This theorem shows that the set of non-solutions cannot escape to infinity. In particular, every element of $\mathcal{I}_n^3$ must lie in a uniformly bounded region of $(0, +\infty)^3$. Consequently, once one excludes some sufficiently large finite box from consideration, the generalised Damascus inequality necessarily holds. 

Now let us present the proofs of these two theorems.

\noindent
{\bf Proof of Theorem \ref{3farfrom1}.} Fix $n\in\mathbf{N}_{\ge 4}$. By symmetry in $(x, y, z)$, it suffices to show that
\[\mathcal{I}_n^3\cap(\mathbf{R}^2\times[1-\delta, 1+\delta])=\emptyset\]
for some $\delta>0$. Throughout this proof, we work with the sup norm $\|\cdot\|_{\infty}$ on $\mathbf{R}^3$. Let
\[\Pi:=\{(x, y, z)\in\mathbf{R}^3 : z=1\}.\]
On $\Pi\cap\mathcal{H}_3$, we have the relation $y=1/x$, so by the case $m=2$ in Section \ref{subsec2.1}, we know that $S_n^3(x, 1/x, 1)\leqslant 0$ for all $x>0$ with equality if and only if $x=1$. Moreover, it follows from the explicit formula that
\[\lim_{x\to 0^+}S_n^3\left(x, \frac{1}{x}, 1\right)=\lim_{x\to+\infty}S_n^3\left(x, \frac{1}{x}, 1\right)=-1.\]
Hence there exist numbers $0<\rho<1/4<4<\rho'$ such that 
\[x\in(0, \rho)\cup(\rho', +\infty)\quad\Rightarrow\quad S_n^3\left(x, \frac{1}{x}, 1\right)<-\frac{7}{8}.\]
We now split the argument into two parts: the ``far" region $x\in(0, \rho)\cup(\rho', +\infty)$ and the ``middle" region $x\in[\rho, \rho']$.

\textit{Scenario 1:} Far. Let us only prove the case $x\in(0, \rho)$ since the case $x\in(\rho', +\infty)$ can be done similarly. First, we note the following facts regarding the function
\[f(x):=\frac{x^n-1}{x^{n+1}+1}.\]
\begin{itemize}
    \item $f$ is increasing and negative on $(0, 1)$.
    \item $f$ is decreasing and positive on $(2, +\infty)$.
\end{itemize}
They can be shown by usual calculus argument, so we will omit the proof. Now fix the reference point $(\rho, 1/\rho, 1)$. By continuity of $S_n^3$ on $\mathbf{R}^3$, there exists $\varepsilon_f\in(0, 1/100)$ such that
\begin{align}\label{yoo}
    S_n^3(x, y, z)<-\frac{3}{4}, \quad \forall (x, y, z)\in B_{\infty}\left(\left(\rho, \frac{1}{\rho}, 1\right), \varepsilon_f\right).
\end{align}
\begin{claim}\label{Claim1}
     $S_n^3$ is negative over $B_{\infty}((\tilde{\rho}, 1/\tilde{\rho}, 1), \varepsilon_f)$ for any $0<\tilde{\rho}<\rho$.
\end{claim}

Analogously, one can also find $\varepsilon_{f'}>0$ such that $S_n^3$ is negative over $B_{\infty}((\tilde{\rho}, 1/\tilde{\rho}, 1), \varepsilon_{f'})$ for any $\tilde{\rho}\geqslant\rho'$.
By proving this claim, we will yield 
\[\mathcal{I}_n^3\cap(((0, \rho)\cup(\rho', +\infty))^2\times[1-\varepsilon_f, 1+\varepsilon_f])=\emptyset,\]
where we may assume without loss of generality that $\varepsilon_f\leqslant\varepsilon_{f'}$. The proof of this claim will be provided later.

\textit{Scenario 2:} Middle. On $\Pi\cap\mathcal{H}_3$ with $x\in[\rho, \rho']$, we have the compact arc
\[\Gamma:=\left\{\left(x,\frac{1}{x}, 1\right) : x\in[\rho, \rho']\right\}.\]
By the strict local maximality of $(1, 1, 1)$ for $S_n^3$ on $\mathcal{H}_3$ (see Section \ref{subsec2.2}), there exists $\varepsilon'>0$ such that
\[S_n^3(x, y, z)<0, \quad\forall (x, y, z)\in B_{\infty}((1, 1, 1), \varepsilon')\cap\mathcal{H}_3, (x, y, z)\neq(1, 1, 1).\]
In particular, we have $S_n^3(x, 1/x, 1)<0$ for each point $(x, 1/x, 1)\in\Gamma\setminus\{(1, 1, 1)\}$. By continuity, for each such $x$, there exists $\varepsilon_x>0$ such that $S_n^3<0$ on $B_{\infty}((x, 1/x, 1), \varepsilon_x)$. Thus the family
\[\left\{B_{\infty}\left(\left(x, \frac{1}{x}, 1\right), \varepsilon_x\right) : \left(x, \frac{1}{x}, 1\right)\in\Gamma\setminus\{(1, 1, 1)\}\right\}\cup\{B_{\infty}((1, 1, 1), \varepsilon')\}\]
forms an open cover of $\Gamma$. Since $\Gamma$ is compact, there exists a finite subcover
\[\left\{B_{\infty}\left(\left(x_k, \frac{1}{x_k}, 1\right), \varepsilon_k\right)\right\}_{k=1}^N\subseteq \left\{B_{\infty}\left(\left(x, \frac{1}{x}, 1\right), \varepsilon_x\right)\right\}\cup\{B_{\infty}((1, 1, 1), \varepsilon')\}\]
such that each ball in this family lies in the region where $S_n^3\leqslant 0$. Let
\[\varepsilon_{m}=\min\{\varepsilon', \varepsilon_1, \ldots, \varepsilon_N\}>0.\]
Then for any $(x, y, z)\in\mathcal{H}_3$ with $|z-1|\leqslant\varepsilon_m$ and $x\in[\rho, \rho']$, the point $(x, y, z)$ lies in one of these finitely many balls. Hence
\[x\in[\rho, \rho'], |z-1|\leqslant\varepsilon_m, xyz=1\quad\Rightarrow\quad S_n^3(x, y, z)\leqslant 0.\]

Finally, define $\delta=\min\{\varepsilon_f, \varepsilon_m\}$. Let $(x, y, z)\in\mathcal{H}_3$ satisfy $|z-1|\leqslant\delta$. Then either $x\in(0, \rho)\cup(\rho', +\infty)$, in which Scenario 1 gives $S_n^3(x, y, z)<0$, or $x\in[\rho, \rho']$, in which Scenario 2 gives $S_n^3(x, y, z)\leqslant 0$. So in all cases, $(x, y, z)\not\in\mathcal{I}_n^3$.\hfill$\square$

\hfill

\noindent
\textbf{Proof of Claim \ref{Claim1}.} Fix $0<\tilde{\rho}<\rho$, and take an arbitrary point $(x, y, z)\in B_{\infty}((\tilde{\rho}, 1/\tilde{\rho}, 1), \varepsilon_f)$. So
\[|x-\tilde{\rho}|\leqslant\varepsilon_f, \quad \left|y-\frac{1}{\tilde{\rho}}\right|\leqslant\varepsilon_f, \quad |z-1|\leqslant\varepsilon_f.\]
Define its ``parallel" point
\[(x', y', z'):=\left(x+\rho-\tilde{\rho}, y+\frac{1}{\rho}-\frac{1}{\tilde{\rho}}, z\right).\]
First, we observe that
\[x'-\rho=x-\tilde{\rho},\quad y'-\frac{1}{\rho}=y-\frac{1}{\tilde{\rho}}, \quad z'-1=z-1.\]
Thus we have
\[\left\|(x', y', z')-\left(\rho, \frac{1}{\rho}, 1\right)\right\|_{\infty}=\left\|(x, y, z)-\left(\tilde{\rho}, \frac{1}{\tilde{\rho}}, 1\right)\right\|_{\infty}\leqslant\varepsilon_f.\]
Hence $(x', y', z')\in B_{\infty}((\rho, 1/\rho, 1), \varepsilon_f)$, and by \eqref{yoo}, we see that
\begin{align*}
    S_n^3(x', y', z')=f(x')+f(y')+f(z')<-\frac{3}{4}.
\end{align*}
Next, we compare $(x, y)$ and $(x', y')$. From
\[|x-\tilde{\rho}|\leqslant\varepsilon_f, \quad 0<\tilde{\rho}<\rho<\frac{1}{4},\quad \varepsilon_f<\frac{1}{100},\]
and $\rho>\tilde{\rho}$,
we get
\[x<x+\rho-\tilde{\rho}=x'<\rho+\varepsilon_f<\frac{1}{2}.\]
Similarly, $x\geqslant\tilde{\rho}-\varepsilon_f>0$. 
So $x, x'\in(0, 1)$ with $x<x'$. Because $f$ is increasing on $(0, 1)$, we have $f(x)<f(x')$. For $y$, from $|y-1/\tilde{\rho}|\leqslant\varepsilon_f$ and $\tilde{\rho}<\rho<1/4$, we have
\[y\geqslant\frac{1}{\tilde{\rho}}-\varepsilon_f>\frac{1}{\rho}-\varepsilon_f>4-\frac{1}{100}>3.\]
Also, $y'=y+(1/\rho)-(1/\tilde{\rho})$ and since $0<\tilde{\rho}<\rho$, we have $1/\tilde{\rho}>1/\rho$. Hence $3<y'<y$. Since $f$ is decreasing and positive on $(2, +\infty)$, we obtain $f(y)<f(y')$. For $z$, we have $|z-1|\leqslant\varepsilon_f$ and $z'=z$. So $f(z)=f(z')$. Combining all these facts together, we deduce that
\begin{align*}
    S_n^3(x, y, z)=f(x)+f(y)+f(z)<f(x')+f(y')+f(z')=S_n^3(x', y', z')<0
\end{align*}
for every $(x, y, z)\in B_{\infty}((\tilde{\rho}, 1/\tilde{\rho}, 1), \varepsilon_f)$. This justifies the claim. \hfill $\square$

\hfill

\noindent
{\bf Proof of Theorem \ref{3bound}.} Fix $n\in\mathbf{N}_{\ge 4}$.  For convenience, define
\[f(x):=\frac{x^n-1}{x^{n+1}+1}\]
for $x>0$. Assume to the contrary that $\mathcal{I}_n^3$ is unbounded. Then for every $c>0$, there exists $(x_c, y_c, z_c)\in\mathcal{I}_n^3$ with $x_c<c$ where we assume without loss of generality that $x_c\leqslant y_c\leqslant z_c$. We record two elementary bounds.
\begin{itemize}
    \item For $0<\tilde{c}<1$, 
\begin{align*}
\frac{(\tilde{c}-\tilde{c}^n)(\tilde{c}+1)}{(\tilde{c}^{n+1}+1)(\tilde{c}^2+1)}\geqslant 0  \quad\Rightarrow\quad  f(\tilde{c})=\frac{\tilde{c}^n-1}{\tilde{c}^{n+1}+1}\leqslant\frac{\tilde{c}-1}{\tilde{c}^{2}+1}.
\end{align*}
    \item For  $\tilde{c}> 1$,
\[f(\tilde{c})=\frac{\tilde{c}^n-1}{\tilde{c}^{n+1}+1}\leqslant\frac{\tilde{c}^n}{\tilde{c}^{n+1}}=\frac{1}{\tilde{c}}.\]
\end{itemize}
In addition, a routine calculus argument shows that the map $\tilde{c}\mapsto (\tilde{c}-1)/(\tilde{c}^2+1)$ is increasing on $(0, 1)$, and the map $\tilde{c}\mapsto 1/\tilde{c}$ is decreasing on $(1, +\infty)$.

By Theorem \ref{3farfrom1}, there exists $\delta\in(0, 1)$ such that for every $(x, y, z)\in\mathcal{I}_n^3$, none of the coordinates lies in $(1-\delta, 1+\delta)$. Since $x_c\leqslant y_c\leqslant z_c$, there are two cases to consider.

\textit{Case 1:} $y_c<1$. Then necessarily $y_c\leqslant 1-\delta$, and from $xyz=1$, we have 
\[z_c=\frac{1}{x_cy_c}\geqslant\frac{1}{x_c}\geqslant\frac{1}{c}>1.\]
Using the bounds above and the monotonicity properties, we obtain
\[f(x_c)\leqslant\frac{x_c-1}{x_c^2+1}\leqslant\frac{c-1}{c^2+1},\]
and since $0<y_c\leqslant 1-\delta<1$,
\[f(y_c)\leqslant\frac{y_c-1}{y_c^2+1}\leqslant\frac{(1-\delta)-1}{(1-\delta)^2+1}.\]
Finally, $z_c>1$ gives $f(z_c)\leqslant 1/z_c\leqslant x_c\leqslant c$. Hence
\begin{equation}\label{S1}
    S_n^3(x_c, y_c, z_c)\leqslant\frac{c-1}{c^2+1}+\frac{(1-\delta)-1}{(1-\delta)^2+1}+c.
\end{equation}
As $c\to 0^+$, the right-hand side of \eqref{S1} tends to
\[\frac{0-1}{0+1}+\frac{(1-\delta)-1}{(1-\delta)^2+1}+0=-1-\frac{\delta}{\delta^2-2\delta+2}=-\frac{\delta^2-\delta+2}{\delta^2-2\delta+2}<0.\]
Therefore, by the definition of limit, there exists $c_1\in(0, 1)$ such that any triple $(x, y, z)$ with $0<x\leqslant c_1$, $0<y\leqslant 1-\delta$, $xyz=1$ satisfies $S_n^3(x, y, z)<0$. In particular, for such $c_1$, no triple in $\mathcal{I}_n^3$ can fall into Case 1 with $x_c\leqslant c_1$.

\textit{Case 2:} $y_c>1$. Then by Theorem \ref{3farfrom1}, we have $y_c\geqslant 1+\delta$. Since $y_c\leqslant z_c$ and $y_cz_c\geqslant 1/c$, we have $z_c\geqslant 1/\sqrt{c}$. Using the bounds above,
\[f(x_c)\leqslant\frac{x_c-1}{x_c^2+1}\leqslant\frac{c-1}{c^2+1}\]
and
\[f(y_c)\leqslant\frac{1}{y_c}\leqslant\frac{1}{1+\delta},\quad f(z_c)\leqslant\frac{1}{z_c}\leqslant\sqrt{c}.\]
Thus
\begin{equation}\label{S2}
    S_n^3(x_c, y_c, z_c)\leqslant\frac{c-1}{c^2+1}+\frac{1}{1+\delta}+\sqrt{c}.
\end{equation}
As $c\to 0^+$, the right-hand side of \eqref{S2} tends to
\[-1+\frac{1}{1+\delta}=\frac{-\delta}{1+\delta}<0.\]
Hence there exists $c_2\in(0, 1)$ such that any triple $(x, y, z)$ with
\[0< x\leqslant c_2,\quad y\geqslant 1+\delta,\quad z\geqslant 1+\delta,\quad xyz=1\]
satisfies $S_n^3(x, y, z)<0$. Particularly, for such $c_2$, no triple in $\mathcal{I}_n^3$ can fall into Case 2 with $x_c\leqslant c_2$.

Finally, choose $c_0=\min\{c_1, c_2\}>0$. By the contradiction assumption, there exists $(x_{c_0}, y_{c_0}, z_{c_0})\in\mathcal{I}_n^3$ with $x_{c_0}\leqslant c_0$. By Theorem \ref{3farfrom1}, this point must be in either Case 1 or Case 2. But for $c=c_0$, both cases lead to $S_n^3(x_{c_0}, y_{c_0}, z_{c_0})<0$, contradicting $(x_{c_0}, y_{c_0}, z_{c_0})\in\mathcal{I}_n^3$. Hence we are done. \hfill $\square$

\begin{example}
Figure~\ref{example} illustrates the geometry of the violation set for the case $(m,n)=(3,6)$. The red region represents the non-positive pre-image of $S_6^3(x,y,z)$, separating the region where the inequality holds from where it fails. The blue surface depicts the constraint hypersurface $\mathcal{H}_3$.

If the inequality were valid on the whole constraint for $(m,n)=(3,6)$, the constraint surface would lie entirely in the non-positive region and would therefore be completely hidden behind the red surface. However, three visible ``fat triangular'' patches of $\mathcal{H}_3$ protrude through the red surface. These correspond exactly to triples $(x,y,z)$ with $xyz=1$ for which $S_6^3(x,y,z)>0$, and thus provide a geometric visualisation of counterexamples in this case.

This picture also reflects two qualitative features established in Theorem \ref{3farfrom1} and Theorem \ref{3bound}. First, the non-solution set $\mathcal{I}_n^3$ stays at a positive distance from the coordinate planes $x=1$, $y=1$, and $z=1$: in particular, the green vertical plane $x=1$ does not intersect the blue patches where $S_6^3>0$ as shown in the right figure. Second, $\mathcal{I}_n^3$ is bounded, so all visible non-solutions lie in a finite region around $(1,1,1)$ rather than escaping to infinity. 

\begin{figure}
    \centering
    \includegraphics[width=12cm]{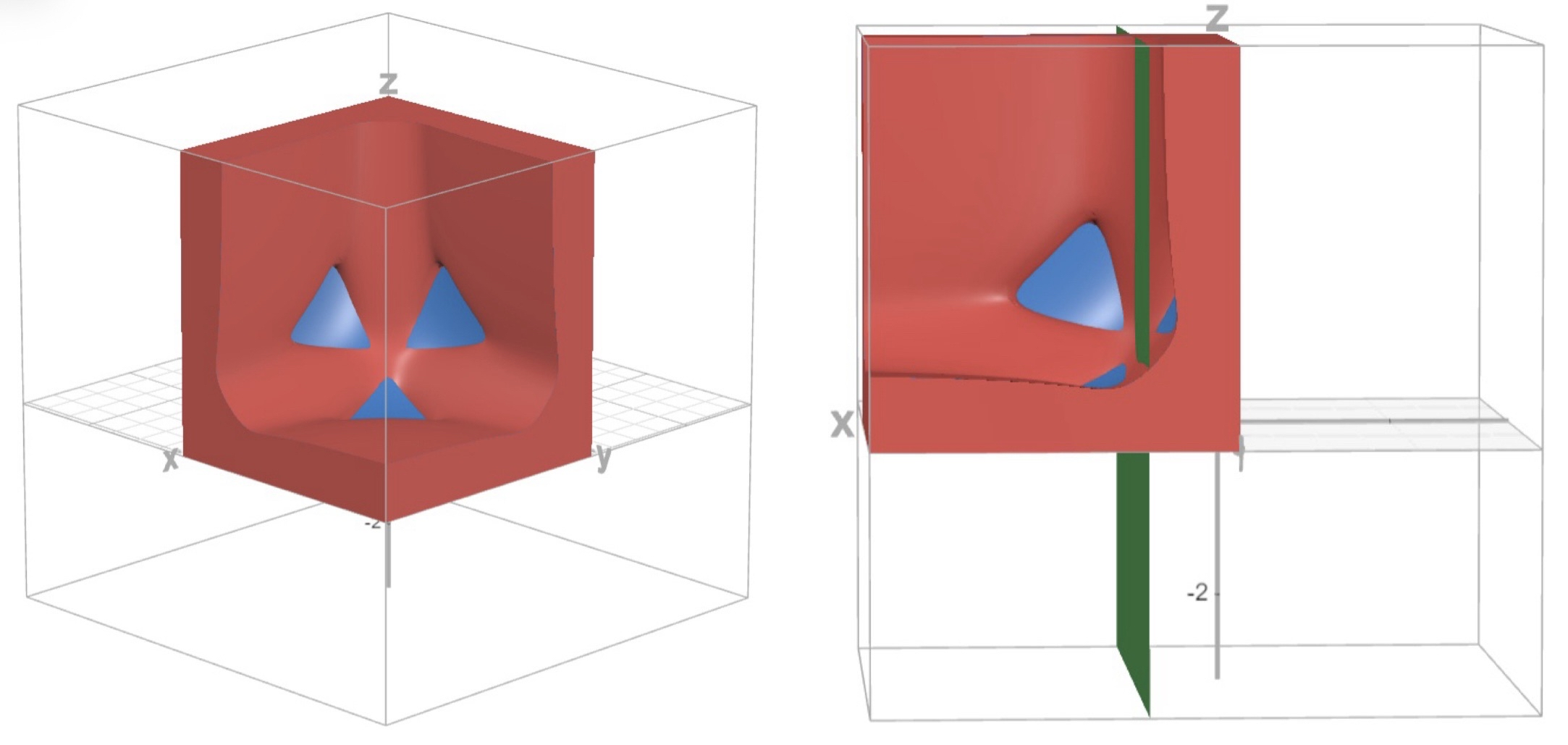}
    \caption{Example for $(m, n)=(3, 6)$}
    \label{example}
\end{figure}
\end{example}


We now investigate how the sets $\mathcal{I}_n^3$ evolve as the exponent
$n$ varies. Motivated by the proofs of Theorem \ref{m=3n>=4}, together with supporting numerical evidence, we propose the following conjecture.

\begin{conjecture}\label{chainofsubset}
For any $n\in\mathbf{N}$, we have $\mathcal{I}_n^3\subseteq\mathcal{I}_{n+1}^3$. 
\end{conjecture}

In other words, this conjecture states that if a triple $(x, y, z)$ fails to satisfy \eqref{varineq} for some $n$, then it also fails for all $n'>n$. This implies that the sequence of sets of non-solutions is nested and increasing. Establishing this conjecture in complete generality appears to be highly nontrivial; nevertheless, if true, it would offer a more refined picture of the structure and progression of non-solutions to the generalised Damascus inequality. We therefore leave this question for future investigation.

In addition to the above, it is natural to ask whether analogous structural properties hold for $\mathcal{I}_n^m$ when $m$ is arbitrary. Can one prove or disprove that 
$\mathcal{I}_n^m$ (i) stays uniformly away from $1$, (ii) is bounded, or (iii) forms an increasing family in $n$? Indeed, the answer to question (i) is negative according to Remark \ref{rmk1}.
We additionally note that the proof of Theorem \ref{3farfrom1} relies crucially on the fact that \eqref{varineq} holds for every $n\in\mathbf{N}$ when $m=2$, a feature that fails for other values of $m$, as shown in the previous sections. 
Because of this, we may likely think that the answer to question (ii) is negative as well. However, we do not know whether the converse ``Theorem \ref{3bound} $\Rightarrow$ Theorem \ref{3farfrom1}" is true for any $m$. Moreover, we have not been able to construct an explicit counterexample yet. If it surprisingly turns out to be true, then it would require a new technique to prove. Thus this question, together with question (iii), remains open. Clarifying these properties of $\mathcal{I}_n^m$ for $m\geqslant 3$ would therefore be an interesting direction for further study, which we leave to the interested readers. 


\section*{Acknowledgements}

We would like to thank the anonymous referee for several constructive comments. The first author also gratefully acknowledges support from the Philip K. H. Wong Foundations Scholarship at the University of Hong Kong.

\hfill

\noindent
\textsc{Chanatip Sujsuntinukul, \\
Department of Mathematics,\\ The University of Hong Kong,\\ Pokfulam, Hong Kong}\\
\textit{E-mail address:} \url{chanatip@connect.hku.hk}

\hfill

\noindent
\textsc{Christophe Chesneau, \\
Department of Mathematics, LMNO,\\
University of Caen-Normandie,\\
14032 Caen, France}\\
\textit{E-mail address:} \url{christophe.chesneau@gmail.com}

\end{document}